\documentclass[11pt]{amsart}
\usepackage{amssymb}
\usepackage{amsmath}
\usepackage{amsrefs}
\usepackage[colorlinks=true]{hyperref}
\usepackage{cleveref}
\usepackage{amscd}
\usepackage{amsthm}
\usepackage{mathtools}
\usepackage{tikz-cd}
\usepackage{bm}
\usepackage{mathrsfs}
\usepackage{enumerate}
\usepackage[margin=2.93cm]{geometry}

\usepackage{xcolor}
\hypersetup{
    colorlinks,
    linkcolor={red!50!black},
    citecolor={green!50!black},
    urlcolor={blue!80!black}
}
\allowdisplaybreaks

\newtheorem{theorem}{Theorem}[section]
\newtheorem{lemma}[theorem]{Lemma}
\newtheorem{proposition}[theorem]{Proposition}
\newtheorem{corollary}[theorem]{Corollary}

\newtheorem{conjecture}[theorem]{Conjecture}
\newtheorem{question}[theorem]{Question}
\newtheorem{fact}[theorem]{Fact}

\newtheorem*{assumption-no-number}{Assumption}
\newtheorem*{corollary*}{Corollary}

\theoremstyle{definition}
\newtheorem{definition}[theorem]{Definition}
\newtheorem{example}[theorem]{Example}

\theoremstyle{remark}
\newtheorem{remark}[theorem]{Remark}

\numberwithin{equation}{section}

\newcommand{\R}{\mathbb{R}}

\newcommand{\N}{\mathbb{N}}
\newcommand{\Z}{\mathbb{Z}}

\renewcommand{\L}{\mathcal{L}}

\newcommand{\tms}{\times}

\newcommand{\rmk}{\begin{remark}}
\newcommand{\ermk}{\end{remark}}
\newcommand{\cor}{\begin{corollary}}
\newcommand{\ecor}{\end{corollary}}
\newcommand{\eq}{\begin{equation}}
\newcommand{\eeq}{\end{equation}}
\newcommand{\eqs}{\begin{equation*}}
\newcommand{\eeqs}{\end{equation*}}
\newcommand{\prop}{\begin{proposition}}
\newcommand{\eprop}{\end{proposition}}
\newcommand{\thm}{\begin{theorem}}
\newcommand{\ethm}{\end{theorem}}
\newcommand{\conj}{\begin{conjecture}}
\newcommand{\econj}{\end{conjecture}}
\newcommand{\lem}{\begin{lemma}}
\newcommand{\elem}{\end{lemma}}
\newcommand{\defi}{\begin{definition}}
\newcommand{\edefi}{\end{definition}}
\newcommand{\ex}{\begin{example}}
\newcommand{\eex}{\end{example}}
\newcommand{\alis}{\begin{align*}}
\newcommand{\ealis}{\end{align*}}
\newcommand{\pf}{\begin{proof}}
\newcommand{\epf}{\end{proof}}
\newcommand{\ali}{\begin{align}}
\newcommand{\eali}{\end{align}}
\newcommand{\qus}{\begin{question}}
\newcommand{\equs}{\end{question}}

\newcommand{\mc}{\mathcal}
\renewcommand{\bf}{\textbf}

\newcommand{\C}{\mathbb{C}}

\newcommand{\sub}{\subseteq}

\newcommand{\ov}{\overline}

\newcommand{\bb}{\mathbb}

\newcommand{\op}{\operatorname}

\renewcommand{\a}{\alpha}

\renewcommand{\d}{\partial}
\newcommand{\e}{\epsilon}

\newcommand{\g}{\gamma}

\newcommand{\s}{\sigma}

\renewcommand{\l}{\lambda}

\newcommand{\fk}{\frak}

\renewcommand{\L}{\Lambda}
\renewcommand{\O}{\Omega}

\renewcommand{\S}{\Sigma}

\newcommand{\Q}{\mathbb{Q}}

\newcommand{\bu}{\bullet}

\renewcommand{\ov}{\overline}

\numberwithin{equation}{section}
\numberwithin{figure}{section}
\title{Some constructions of Weinstein manifolds with chaotic Reeb dynamics}
\author{Laurent C\^{o}t\'{e}}

\begin{document}

\date{\today}  
\begin{abstract} 
The goal of this note is to describe some constructions of Weinstein manifolds with chaotic Reeb dynamics, and to explain how this property can sometimes be detected directly from the skeleton.
\end{abstract}

\maketitle

\section{Introduction}

\subsection{Some problems and questions}\label{subsection:intro-questions}
Let $X$ be a Liouville manifold. For the purpose of studying Reeb dynamics on $\d_\infty X$, it is often useful to consider the \emph{growth rate} of Floer theoretic invariants (such as symplectic cohomology, wrapped Floer cohomology, contact homology\dots) as a function of the action. This growth rate a priori depends on some choices, but one can show that it is well-defined up to a natural notion of scaling equivalence. 

When the growth rate is exponential, this tends to force the Reeb dynamics to be chaotic. This heuristic is made precise e.g.\ in work of Macarini--Schlenk \cite{macarini-schlenk}, Alves \cite{alves}, Alves--Meiwes \cite{alves-meiwes} and others.  McLean \cites{mclean-gafa, mclean} (building on work of Seidel \cite{seidel-biased-view}) also showed that $SH(X)$ (and $HW(K, L)$) cannot grow super-polynomially if $X$ is symplectomorphic to an affine variety; this ultimately comes from the fact that affine varieties can be shown to admit a nice contact form at infinity having tame Reeb dynamics.

These phenomena are now well-understood. Nevertheless, various questions remain which are (to the author's taste) rather natural. Here are some closely related problems and questions which will be taken up in this note.  


\begin{enumerate}
\item\label{item:skeleton-weinstein} Describe conditions on the skeleton of a Weinstein manifold $X$ which guarantee that $X$ is not symplectomorphic to an affine variety.
\item\label{item:subcritical} Suppose that $X$ is obtained by attaching a \emph{critical} handle to $T^*M$, where $M$ is a closed hyperbolic manifold. Describe conditions on the attaching sphere which ensure that $X$ is not symplectomorphic to an affine variety.
\item\label{item:no-exact-lag} Construct a Weinstein manifold which is not symplectomorphic to an affine variety and which does not contain a closed exact Lagrangian.
\item\label{item:exotic-s5} Does there exist a contact structure $\xi_{\op{exotic}}$ on $S^5$ with the property that for every contact form $\op{ker} \a=\xi_{\op{exotic}}$, the associated Reeb flow has positive topological entropy?
\item\label{item:entropy-monodromy} Let $f: X \to \C$ be a Lefschetz fibration on a Weinstein manifold. Consider the partially wrapped Fukaya category $\mc{W}(X, f^{-1}(\infty))$. This is a smooth and proper category whose Serre functor $\mathbf{W}[-1]$ was described by Sylvan \cite{sylvan}. Describe conditions on $X$ (independent of $f$) which guarantee that the Serre functor has positive categorical entropy. 

\end{enumerate}

\subsection{What is known?}

Concerning \eqref{item:skeleton-weinstein}, it is known that $\op{skel}(X, \l)= M$ implies that $(X, \l)= (T^*M, \l_{\op{std}})$. When $M$ is a hyperbolic (or rationally hyperbolic) manifold, it follows from work of McLean \cite{mclean-gafa} that $X$ is not symplectomorphic to an affine variety. On the opposite end of the spectrum, it is not difficult to show that \emph{any closed manifold} of real dimension $n$ can be embedded in the skeleton of a smooth affine variety of complex dimension $n$ (see \Cref{proposition:nash}). To the best of our knowledge, nothing is known beyond these two extreme cases.

Concerning \eqref{item:subcritical} and \eqref{item:no-exact-lag}, McLean \cite{mclean-gafa} has constructed many examples of Weinstein manifolds which are not symplectomorphic to affine varieties. To the best of our knowledge, all of these examples are built as subcritical handle attachments on cotangent bundles (or products of these). 


Concerning \eqref{item:exotic-s5}, Alves and Meiwes \cite{alves-meiwes} proved the analogous result on $S^{2n+1}$ for all $n>2$. However, the case $n=2$ remained open. (The analogous statement when $n=1$ is expected to be false.)

Concerning \eqref{item:entropy-monodromy}, it was observed in \cite{cote-kartal} (see Ex.\ 8.8 and Sec.\ 8.2 in \textit{loc. cit.}) that this property holds for cotangent bundles of hyperbolic and rationally hyperbolic manifolds. 

\subsection{The contribution of this note}
The goal of this note is to study the above questions. The main reason why we are able to say anything new is mostly thanks to the recent work of Ganatra--Pardon--Shende \cite{gps2} on the foundations of wrapped Fukaya categories, which we combine with ideas of other authors including McLean \cite{mclean} and Alves--Meiwes \cite{alves-meiwes}. 


As mentioned above, the standard approach for studying questions of this type is to consider an action filtration on either symplectic cohomology or wrapped Floer cohomology.  One then studies the growth of the hom-spaces with respect to this filtration. When $X=T^*M$, the action filtrations on $SH/HW$ match up with the action filtrations on $\mc{L}M/ \O M$ induced by a choice of metric on $M$, which allows one to perform many computations. However, to go beyond cotangent bundles, one needs to understand how the action filtration behaves under various operations. This is actually quite hard to do. For symplectic cohomology, McLean \cite{mclean-unpublished} described the effect of subcritical handle attachments and products on the filtered growth, and this was already non-trivial. Trying to redo the entire package of \cite{gps2} in the filtered setting seems like a rather daunting undertaking. 

Instead of doing this, we use a basic but important observation due to Alves and Meiwes \cite{alves-meiwes}, namely that one can bound the action-filtered growth of wrapped Floer cohomology from below in terms of the internal multiplicative growth. This is a direct consequence of the fact that the action filtration is sublinear with respect to multiplication. However, the multiplicative growth is a purely algebraic invariant which does not need a filtration, so it is very well behaved under the package of operations in \cite{gps2}.

\subsection{Context and acknowledgements}
This note replaces a previous paper which was roughly three times longer, so let us explain what changed. In the previous version, the aim was to define a categorical notion of ``algebraic hyperbolicity" which could be used to detect chaotic Reeb dynamics and would be well-behaved under various natural operations. However, we subsequently realized that essentially all of the applications can be obtained (and in some cases improved) by much simpler arguments. We also found various errors in the previous version, some of which seem serious. Since the categorical approach appears to be a distraction, we decided to abandon it entirely. 

The author has benefited from conversations, suggestions and correspondences with many mathematicians including Dani \'{A}lvarez-Gavela, Marcelo Alves, Denis Auroux, Georgios Dimitroglou Rizell, Yasha Eliashberg, Yves F\'{e}lix, Kathryn Hess, Helmut Hofer, Yusuf Bar\i\c{s} Kartal, Thomas Massoni, John Pardon, Semon Rezchikov, Lisa Sauermann, Vivek Shende, Kyler Siegel, Zack Sylvan and Umut Varolgunes. The author is particularly grateful to Marcelo Alves for explaining to him the arguments of \Cref{subsection:exotics5}. Part of this research was carried out while the author was a graduate student at Stanford University and supported by a Benchmark Graduate Fellowship. This material is based upon work supported by the National Science Foundation under Grant No. DMS-1926686 (while the author was at the Institute for Advanced Study). 




\section{Main definitions}

\subsection{Topological and categorical entropy}

The topological entropy is one of the simplest measures of complexity of a dynamical system. Intuitively, a system with positive topological entropy can be viewed as chaotic. 

Let $(X,d)$ be a compact metric space. Let $\Phi= \{\phi_t\}$ be a one-parameter family of homeomorphisms and consider the family of metrics \eq d^{\Phi}_T(x,y) = \max_{0 \leq t \leq T} d(\phi_t x, \phi_t y)\eeq depending on $T>0$. Let $B_{\Phi}(x, \e, T)= \{ y \in X \mid d^{\Phi}_T(x,y)\}< \e$ be the ball of radius $\e$ centered at $x$.  Now let $S_d(\Phi, \e, T)  \in \N$ be defined as the smallest number of balls of radius $\e$ needed to cover $X$. Set $h_d(\Phi, \e):= \op{limsup}_{T \to \infty} \frac{1}{T} \op{log} S_d(\Phi, \e, T)$.

The \emph{topological entropy} of the flow $\Phi= \{\phi_t\}$ on the compact metric space $(X, d)$ is \eq h(\Phi)= \lim_{\e \to 0} h_d(\Phi, \e).\eeq It can be shown \cite[Sec.\ 3.1]{k-h} that this quantity only depends on $\Phi$ and the topology induced by $d$. 

\ex 
It follows from the definition that $h(\Phi)=0$ if $\phi_t$ is an isometry for all $t$, or if it is periodic. More surprisingly, it can be shown \cite[Sec.\ 3.3]{k-h} that the topological entropy of any gradient flow on a compact manifold vanishes.  In contrast, the geodesic flow on a connected, orientable surface $\S$ endowed with a metric of constant curvature has positive topological entropy if $g \geq 2$ and vanishing topological entropy otherwise \cite[Sec.\ 3.2.5]{paternain}. 
\eex

There is a categorical analog of the topological entropy which has introduced by Dimitrov--Haiden--Kontsevich--Katzarkov \cite{d-h-k-k}. Given a smooth an proper $A_\infty$ category $\mc{C}$ and an endofunctor $F$, the \emph{(categorical) entropy} of $F$ is the quantity:
\eq h(F):= \lim_{n \to \infty} \log \op{dim} \hom_{H(\mc{C})}(G, F^n G)  \in [-\infty, \infty) \eeq
where $G$ is any (split)-generator of $\mc{C}$. This quantity can be shown to be independent of the choice of $G$. (More generally, the categorical entropy can be defined for any endofunctor acting on a triangulated category.)

\subsection{Liouville manifolds}

We assume that the reader is familiear with the definition of Liouville manifolds/sectors and their wrapped Fukaya category; see \cite{gps1} for the relevant background. Unless otherwise indicated, we will always consider the wrapped Fukaya category with $\Z/2$-gradings and $\Q$-coefficients. In particular, all Liouville manifolds/sectors and all Lagrangians in the Fukaya category are equipped with appropriate orientation/grading data. (The reason for taking $\Q$-coefficients is that we sometimes want to appeal to computations from rational homotopy theory.)

\subsection{Multiplicative growth}

Let $A$ be an associative $k$-algebra. Given a finite set of morphisms $\S \sub \op{Mor}(\mc{C})$ and an integer $n \geq 1$, let $W_{\S}(n)$ be the span of all composable words of length at most $n$ in the elements of $\S$, that is 
\eq W_\S(n)= \op{span}_k\{a_1\dots a_l \mid a_i \in \S, l \leq n \} \sub A.\eeq


\defi \label{definition:alg-growth}
An associative algebra $A$ is said to have \emph{exponential (multiplicative) growth} if 
\eq \op{limsup}_n \frac{1}{n} \log \dim_k W_\S(n) >0 \eeq
for some finite set of elements $\S \sub A$. An object $K$ in an $A_{\infty}$ category $\mc{C}$ is said to have exponential algebraic growth if $H^\bu(\mc{C})(K, K)$ has exponential growth.
\edefi

\ex 
A group $G$ has exponential growth (in the usual sense of e.g.\ \cite[Sec.\ 1]{grigorchuk-pak}) iff the group ring $\Z[G]$ has exponential multiplicative growth. 
\eex

We now introduce the following terminology:
\defi
A Liouville sector is said to be \emph{big} is there exists a Lagrangian disk $K \in \mc{W}(X)$ (equipped with appropriate orientation/grading data) such that $HW^\bu(K,K)$ has exponential multiplicative growth. 
\edefi

The terminology ``big" is intended to be purely auxiliary (i.e.\ the reason for introducing this adjective is to make certain statements in this paper less wordy, and we trust that it will not cause confusion with unrelated notions of big-ness in other areas of mathematics.)

The easiest examples of big Liouville sectors are the following cotangent bundles.

\ex\label{example:hyperbolics}
Let $M$ be an aspherical compact manifold with boundary. If $\pi_1(M)$ has exponential growth, then $T^*M$ is big.  Examples include: any closed manifold admitting a metric with all sectional curvatures strictly less than zero (see \cite[Thm.\ 2]{milnor}). 
\eex

\ex\label{example:q-hyperbolics}
Let $M$ be a (say simply-connected, spin) manifold and let $F \sub T^*M$ be a cotangent fiber. Then $CW^\bu(F, F) \simeq C_{-\bu}(\O M;\Q)$. Hence $T^*M$ is big if $H_\bu(\O M; \Q)$ has exponential multiplicative growth. Morally, one expects this property to hold whenever $M$ is rationally hyperbolic; this would in particular follow from a positive solution to the celebrated Avramov--F\'{e}lix conjecture \cite[Sec.\ 39.4]{q-homotopy}. Here are some additional assumptions on $M$ which ensure that $H_\bu(\O M; \Q)$ has exponential multiplicative growth (and hence $T^*M$ is big): 
\begin{itemize}
\item[(i)] $M$ has the rational homotopy type of a wedge of $m \geq 2$ spheres (see \cite[Sec.\ 33(c); Ex.\ 1]{q-homotopy}); 
\item[(ii)] $H_{2i}(M;k)=0$ for all $i >0$ and $M$ does not have the rational homotopy type of a sphere (see \cite[Sec.\ 24(c)]{q-homotopy}); 
\item[(iii)] $M$ is a connected sum of spin manifolds, each of whose rational cohomology algebra requires at least two generators (see \cite[Cor.\ 3.4 and Rmk.\ 3.5]{f-o-t});
\item[(iv)] $M$ is rationally hyperbolic and $H_\bu(\O X; k)$ is finitely generated as an algebra. 
\end{itemize}
\eex

\subsection{Consequences of exponential multiplicative growth}
The following theorem is a restatement of work of McLean \cite{mclean}, Alves--Meiwes \cite{alves-meiwes} and the author and Kartal \cite{cote-kartal}. 
\thm
Suppose that $(X, \l)$ is big. Then:
\begin{enumerate}
\item\label{item:affine-variety-symp} $X$ is not symplectomorphic to an affine variety;
\item\label{item:link-cont} $\d_\infty X$ is not contactomorphic to the link of an isolated complex singularity;
\item\label{item:am-result} for any contact form on $(\d_\infty X, \op{ker} \l)$, the Reeb flow has positive topological entropy;
\item\label{item:cat-entropy} if $X$ is Weinstein and $\fk{f} \sub \d_\infty X$ is the page of a Weinstein open book, then $\mathbf{W}^{-1}, \mathbf{M}^{-1}$ have positive categorical entropy. 
\end{enumerate}
\ethm
The functors $\mathbf{W}^{-1}: \op{Tw} \mc{W}(X, \fk{f})$ and $\mathbf{M}^{-1}: \op{Tw} \mc{W}(\hat{\fk{f}}) \to \op{Tw} \mc{W}(\hat{\fk{f}})$ were described by Sylvan \cite{sylvan} and are, respectively, the ``wrap-once negatively" and ``inverse monodromy" auto-equivalences. Modulo some folkloric compatibility statements, $\mathbf{W}^{-1}$ is the Serre functor of the Fukaya--Seidel category $\mc{W}(X, \fk{f})$ (as a consequence of $\mathbf{W}^{-1}$ being a Serre functor, $\mathbf{W}$ has positive entropy iff $\mathbf{W}^{-1}$ does; see the discussion in \cite[Sec.\ 8.2]{cote-kartal}). 

\pf
\eqref{item:am-result} was proved by Alves--Meiwes in \cite[Thm.\ 1.7]{alves-meiwes}. Let us consider \eqref{item:cat-entropy}. First, we endow $\mc{W}(X)$ with an integral action filtration according to \cite[Def.\ 5.1]{cote-kartal}. (The idea is to present the wrapped Fukaya category as the localization of some (homologically smooth and proper) partially wrapped Fukaya category, and to use the filtration induced by this presentation). Concretely, this means that for every $K, L \in \mc{W}(X)$, there is a filtration $F^pHW^\bu(K, L)$ which is subadditive with respect to multiplication. The assumption that $X$ is big immediately implies that $F^pHW^\bu(K, K)$ grows exponentially for some Lagrangian disk $K$ (by subadditivity). This implies \eqref{item:cat-entropy} by \cite[(8.9), (8.9)]{cote-kartal} (note that the growth of $F^pHW^\bu(K, K)$ is denoted by $\gamma_{K,L}$ in \cite{cote-kartal}). 

Let $H$ be a (say time-independent) Hamiltonian which is cylindrical at infinity with positive slope and let $\phi_{nH}$ denote the flow of $nH$ for $n \in \N$. Then it is proved in \cite[Thm.\ 7.2]{cote-kartal} that the growth of the sequence $\gamma^{ham}_{K,L}:= \op{dim}(HF^\bu(\phi_{nH}K, L) \to HW^\bu(K, L))$ agrees (up to shift and scaling) with the growth of $F^pHW^\bu(K, L)$. Now $\gamma^{ham}_{K,L}$ is precisely the quantity considered by \cite[Def.\ 2.4]{mclean}, so \eqref{item:affine-variety-symp} and \eqref{item:link-cont} follow from \cite[Thm.\ 1.2]{mclean}. 
\epf

\subsection{Miscellaneous facts}

The following proposition follows from \cite[Thm.\ 3.27]{AGEN3} (note however that this particular statement certainly does not need \cite{AGEN3}; it is for instance already contained in Sec.\ 12.3 in \cite{c-e}). 
\prop \label{proposition:abundance}
Let $(X,\l)$ be a Weinstein manifold. Let $M$ be a compact manifold (possibly with boundary) and suppose that $i: M \to \op{Skel}(X,\l)$ is a smooth embedding into the Lagrangian stratum. Then after possibly deforming $\l$ in a compact set, $i$ is covered by an inclusion of Liouville sectors $\tilde{i} (T^*M, \l_{\op{std}}) \hookrightarrow (X,\l)$. 
\qed
\eprop

In particular, it follows from \Cref{proposition:abundance} that any point in the Lagrangian stratum is contained the image of some sectorial embedding $(T^*M, \l_{\op{std}}) \to (X,\l)$ after possibly deforming $\l$. Note that the proof of \Cref{proposition:abundance} uses the assumption that $(X, \l)$ is Weinstein in an essential way. It is unknown whether a similar result holds for arbitrary Liouville manifolds.

The next proposition, which was promised earlier, shows that \emph{any} closed $n$-manifold embeds into the skeleton of a smooth affine variety of complex dimension $n$. 
\prop \label{proposition:nash}
Let $M$ be a closed manifold of real dimension $n$. For any $N \geq 2n+1$, there exists a smooth affine variety $X \sub \C^N$ of complex dimension $n$ with the following properties:
\begin{itemize}
\item $(X, \l|_X)$ is a Weinstein manifold, where $\l=\sum_{i=1}^N (x_i dy_i - y_i dx_i)$ is the standard radial Liouville form; 
\item $M$ is diffeomorphic to a connected component of the skeleton $\op{Skel}(X, \l|_X)$.
\end{itemize}
\eprop

\pf
According to a classical theorem of Nash (see \cite[Thms.\ 2 and 3]{kollar}), there exists a collection $f_1,\dots,f_k$ of real polynomials in $N$ variables such that $M$ is a connected component of the zero set $V_{\R}(f_1,\dots,f_k) \sub \R^N$. Viewing the $f_i$ as complex polynomials, we set $X:= V_{\C}(f_1,\dots,f_n) \sub \C^n$. After possibly perturbing the $f_i$, we may assume that $X$ is a smooth complex algebraic variety and that the function $\phi(z_1,\dots,z_n) = |z_1|^2+\dots + |z_n|^2$ is Morse and exhausting. One can check that $\l= d^\C\phi$, from which it follows that $(X, \l|_X, \phi)$ is a Weinstein structure (see \cite[Sec.\ 1.1]{c-e}).

Observe that $X$ is stable under complex conjugation, which sends $\l \mapsto - \l$ and $d\l \mapsto - d\l$. It follows that the Liouville vector field $v_X$ on $(X, \l|_X)$ is preserved by complex conjugation. In particular, $v_X$ is tangent to the real locus $X_{\R} \sub X$. Since $M$ is a connected component of $X_\R$, it follows immediately that $M$ is contained in the skeleton of $(X,\l|_X)$. It is similarly straightforward to check that $X$ has complex dimension $n$: indeed, since $T_mX$ is also fixed by complex-conjugation for any $x \in M$, it follows that $\op{dim}T_mX= 2 \op{dim} T_mM$.
\epf

\section{Operations which preserve multiplicative growth}

\subsection{Subcritical and flexible handle attachments}

Given a Liouville sector $(X,\l)$ of dimension $2n$, recall that a Liouville hypersurface $F \sub \d X$ is said to be \emph{subcritical} if its its skeleton is contained in the smooth image of a second countable manifold of dimension $<n-1$.

\prop 
\label{proposition:subcritical-embedding}
Let $(X_1,\l_1)$ and $(X_2,\l_2)$ be Liouville manifolds. Suppose that $(X, \l)$ is a Liouville manifold obtained by gluing $X_1, X_2$ along a subcritical Weinstein hypersurface $F$. Then 
\eq \mc{W}(X_1) \sqcup \mc{W}(X_2) \to \mc{W}(X) \eeq
is a pre-triangulated equivalence. Hence $\mc{W}(X)$ is big if $\mc{W}(X_1)$ or $\mc{W}(X_2)$ is big. 
\eprop

\pf
Since $F$ is subcritical we have $\mc{W}(F)=0$. We also have $\mc{W}(X_i, \fk{c}_F) = \mc{W}(X_i)$ by stop removal (see \cite[Thm.\ 1.16]{gps2}). Now apply \cite[Thm.\ 1.20]{gps2}.
\epf

\lem \label{lemma:acyclic-loose-linking}
Let $(X, \l)$ be a Weinstein manifold and let $F \sub \d_\infty X$ be a Liouville hypesurface. Suppose that $\op{Skel}(F) \sub \d_\infty X$ is a (smooth) sphere and let $D \sub \mc{W}(X, \L)$ be the linking disk. If $\L$ is a loose Legendrian sphere, then $CW^\bu(D, D)_{X, \L}$ is acyclic. 
\elem

\pf
This following from work of Ekholm--Lekili \cite{ekholm-lekili} who prove (see Thm.\ 2 and the surrounding discussion) that $CW^\bu(D, D)_{X, \L}$ is isomorphic to the Chekanov--Eliashberg dg algebra of $\L$ with loop-space coefficients, which is acyclic since $\L$ is loose. (An alternative argument, which has the advantage of remaining entirely in the framework of \cite{gps2}, will appear in forthcoming work of O.\ Lazarev.)
\epf

\prop \label{proposition:flexible-faithful}
Let $(X_1, \l_1)$ and $(X_2, \l_2)$ be Liouville manifolds. Suppose that $(X,\l)$ is obtained by gluing $X_1, X_2$ along the a Weinstein hypersurface $F$. If $\op{Skel}(F) \sub \d X_1$ is a loose Legendrian sphere $\L$ and $(X_1, \l_1)$ is Weinstein up to deformation, then there exists a fully faithful embedding $\mc{W}(X_1) \to \mc{W}(X)$. Hence $\mc{W}(X)$ is big if $\mc{W}(X_1)$ is big. 
\eprop

\pf
Since $F$ is Weinstein, $F= T^*\L$ up to deformation. The map $\mc{W}(F) \to \mc{W}(\check{X}_1) = \mc{W}(X_1, F)= \mc{W}(X_1, \L)$ sends the cotangent fiber to the linking disk of $\L$. In particular $\mc{W}(F) \to \mc{W}(\check{X}_1)$ has acyclic image according to \Cref{lemma:acyclic-loose-linking}. The proposition follows. 
\epf

\subsection{Products}

\prop\label{proposition:products}
Let $X_0, X_1$ be Weinstein manifolds and let $X=X_0 \tms X_1$. Suppose that the following two conditions hold:
\begin{itemize}
\item $X_0$ or $X_1$ is big;
\item $\mc{W}(X_0)$ and $\mc{W}(X_1)$ are both nonzero (or equivalently, $SH(X_0), SH(X_1)$ are nonzero by \cite{ganatra}). 
\end{itemize}
Then $X$ is big.
\eprop

\pf
According to \cite[Thm.\ 1.5]{gps2}, there is a fully faithful embedding
\eq\label{equation:kunneth-embedding} \mc{W}(X) \otimes \mc{W}(X') \hookrightarrow \mc{W}(X \tms X'). \eeq 
If $X_0$ is big, then there exists a disk $K \in \mc{W}(X)$ with exponential multiplicative growth. By assumption (ii), $\mc{W}(X_1)$ is non-zero and Weinstein. Hence, since Weinstein manifolds are generated by cocores (by deep work of Chantraine--Dimitroglou Rizell--Ghiggini--Golovko \cite{cdrgg} and independently Ganatra--Pardon--Shende \cite{gps2}), it follows that there exists a cocore disk $L \in \mc{W}(X_1)$ such that $HW(L_1, L_1) \neq 0$. Now $HW^\bu(K \tms L, K \tms L)$ contains $HW^\bu(K, K) \otimes HW^\bu(L, L)$ as a subalgebra by \eqref{equation:kunneth-embedding}, and hence contains $HW^\bu(K, K)= HW^\bu(K, K) \otimes e_L$ as a subalgebra. This proves the claim. 
\epf

\subsection{Covariant functoriality}

Recall that a Liouvile domain is just the truncation of a Liouville manifold. Similarly, the truncation of a Liouville sector shall be called a (sectorial) \emph{Liouville domain with corners}. More precisely, a (sectorial) Liouville domain with corners is an exact symplectic manifold with corners $(V, \l)$ whose boundary admits a decomposition $\d V= \d_h V \cup \d_v V$. We require that $(\d_h V, \l)$ is a contact manifold with convex boundary. We also require that the Liouville flow preserve $V$ in a neighborhood of $\d_h V$, so that there is a natural decomposition $([-\e, 0] \tms \d_hV)$.

Let $j: X_0^{\op{in}} \hookrightarrow X_0$ be an inclusion of (honest) Liouville domains (i.e.\ a proper embedding which strictly pulls back the Liouville form). Let $X^{\op{in}}$ (resp.\ $X$) denote their completions. Note that $X^{\op{in}}$ does not in general admit a proper embedding into $X$. 

Suppose that $i^{\op{in}}: V \hookrightarrow X_0^{\op{in}}$ and $i: V \hookrightarrow X_0$ are inclusions of Liouville domains with corners with the following properties: (i) $i= j \circ i^{\op{in}}$; (ii) $i^{\op{in}}$ and $i$ induce inclusions of Liouville sectors $\hat{V} \hookrightarrow X^{\op{in}}$ (resp.\ $\hat{V} \hookrightarrow X$). 

\prop\label{proposition:covariant-functoriality}
Under the above hypotheses, the following diagram commutes
\eq
\begin{tikzcd}
& \op{Tw} \mc{W}(X) \ar{d} \\
\mc{W}(\hat{V}) \ar{ur} \ar{r} & \op{Tw} \mc{W}(X^{\op{in}}) 
\end{tikzcd}
\eeq
\eprop
\pf
The wrapped Fukaya category is functorial under inclusions of stopped Liouville sectors. But the Viterbo restriction functor is constructed in \cite[Sec.\ 8.2]{gps2} as a composition of such inclusions. 
\epf

\cor[of \Cref{proposition:covariant-functoriality}]\label{corollary:cov-func}
If there exists a Lagrangian disk $L \in \mc{W}(\hat{V})$ such that the image of $HW_{\mc{W}(\hat{V})}(L, L) \to HW_{\mc{W}(X_0)}(L, L)$ has exponential multiplicative growth, then $(X, \l)$ is big. 
\qed
\ecor

\section{Applications}

\subsection{Recognizing big Weinstein manifolds from their skeleton}

The following proposition gives a criterion to check that a Weinstein manifold is big in terms of the stratified topology of its skeleton.

\prop \label{proposition:plumbings-immersions}
Let $(X, \l)$ be a Weinstein manifold of dimension $2n$. Let $M$ be a closed, connected manifold of dimension $n$, possibly with boundary. Suppose that there exists a smooth embedding $f: M \to X$ which satisfies the following two conditions:
\begin{itemize}
\item[(i)] the image of $f$ is contained in the (open!) Lagrangian stratum of $\op{Skel}(X,\l)$;
\item[(ii)] there is an extension
\eq
\begin{tikzcd}
M \ar[r, hook] \ar[dr, "f" ']&  K \ar[d, "\tilde{f}"]\\
& X 
\end{tikzcd}
\eeq
where $\tilde{f}: K \to X$ is an exact Lagrangian embedding. 
\end{itemize}
If the image of the induced map $H_\bu(\O M; \Q) \to H_\bu(\O K ; \Q)$ has exponential multiplicative growth, then $X$ is big.  
\eprop

\begin{figure}
\includegraphics[width=80mm]{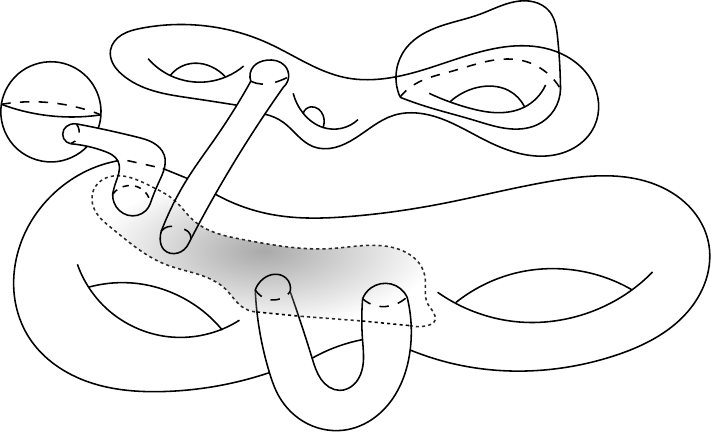}
\vspace{5mm}
\caption{This a priori complicated arboreal skeleton admits an obvious embedding of the genus $2$ surface with the shaded region removed. Its Weinstein thickening  (which exists in light of \cite{AGEN3}) can thus be recognized as big by \Cref{proposition:plumbings-immersions}.}
\label{figure:critical-gluing}
\end{figure}

\pf
Applying \Cref{proposition:abundance}, we may assume (after possibly deforming $\l$ while keeping $\op{skel}(X, \l)$ fixed) that $f$ extends to a \emph{strict} Liouville sector embedding $\s: T^*M \hookrightarrow X$. By the relative version of Moser's theorem (and after possibly slightly shrinking $M$ inwards and further deforming $\l$), $\tilde{f}$ extends to a strict Liouville embedding $\tilde{\s}: X_0:=T^*_{\leq 1} K \hookrightarrow X$ which agrees with $\s$ on their common domain of definition. We are now in the setting of \Cref{proposition:covariant-functoriality} with $V=T^*M$. The result thus follows from \Cref{corollary:cov-func} and the standard isomorphism between the wrapped Floer theory of a cotangent fiber and the homology of the based loopspace (see e.g.\ \cite[Cor. 6.1]{gps3} for the most general version, where it is also shown that this isomorphism is compatible with codimension zero inclusions). 
\epf

%

\rmk
Both assumptions in the statement of \Cref{proposition:plumbings-immersions} are needed. Indeed, \Cref{proposition:nash} implies that any manifold whatsoever can be embedded in the skeleton of an algebraic variety (which is therefore not big): this explains why one must impose that the image of $f$ is contained in the Lagrangian stratum of $\op{Skel}(X,\l)$. To see that (ii) is also needed, let $N$ be a closed hyperbolic manifold of dimension $>2$ and let $g$ be a Morse function with a single maximum. Consider the induced handle decomposition on $X= T^*N$. Now let $X'$ be obtained by ``loosifying" the top handle attachment, that is, by attaching the top handle along a stabilization of the original Legendrian knot. Then the skeleton of $X'$ admits an embedding of $N-B$ into the Lagrangian stratum. However, $X'$ is a flexible Weinstein manifold, so its wrapped Fukaya category is zero. In particular, it is not big.
\ermk


\ex\label{example:skeleton-hyp}
Suppose $M$ is a compact hyperbolic manifold (more generally, $M$ admits a metric with all sectional curvatures strictly less than zero). Then $H^\bu(\O_* (M-D)) \to H^\bu(\O_*M)= k[\pi_1(M)]$ is surjective. Now $k[\pi_1(M)]$ has exponential multiplicative growth by \Cref{example:hyperbolics}. Hence $X$ is big. 
\eex

\ex\label{example:skeleton-rat-hyp}
Suppose $M$ is a rationally hyperbolic compact spin manifold whose rational cohomology algebra requires at least $2$ generators, and such that $H_*(\O M;\Q)$ has exponential multiplicative growth.  (Morally, ``most" spin manifolds satisfy this condition; to get an explicit example, take (i) or (iii) in \Cref{example:q-hyperbolics}). 

According to \cite[Thm.\ 3.3]{f-o-t}, the inclusion $M-D \to M$ induces a surjection $\pi_*(\O(M-D)) \otimes \Q \to \pi_*(\O M) \otimes \Q$. Passing to the universal enveloping algebra, we get a surjection of associative algebras $H_*(\O(M-D); \Q) \to H_*(\O M;\Q)$ (the universal enveloping algebra functor is right exact). Hence $X$ is big. 
\eex


\ex[Critical handle attachments]
\Cref{proposition:plumbings-immersions} immediately allows one to detect critical handle attachments which do not destroy ``bigness". For example, suppose that $M$ is hyperbolic and let $X$ be obtained by attaching a critical handle to $T^*M$. If the front projection of the attaching sphere lands inside a disk, then it follows from \Cref{example:skeleton-hyp} that $X$ is big.  (More generally, the same conclusion holds if the front projection is contained in an open set $U \sub M$ such that $\pi_1(M-U)$ has exponential growth). 
\eex

\subsection{Plumbings and immersions}

As an illustration of \Cref{proposition:plumbings-immersions}, we consider skeleta with ``simple normal crossing" singularities.  It will be convenient to let $\mathscr{E}$ denote the set of all closed manifolds which satisfy the assumptions of \Cref{example:skeleton-hyp} or \Cref{example:skeleton-rat-hyp}. 

\prop \label{proposition:plumbings} 
Suppose that $(X,\l)$ is Weinstein and that $\op{Skel}(X,\l)$ is a union of closed Lagrangian submanifolds $L_1,\dots,L_l$ with the property that $L_i \bigcap (\bigcup_{j \neq i} L_i) \sub L_i$ is contained in a ball. If $L_i \in \mathscr{E}$ for some $i$, then $(X,\l)$ is big. 
\eprop

\pf
Up to reordering the $L_i$, we may assume that $L_1 \in \mathscr{E}$. By assumption, there exists a ball $B \sub L_1$ such that $B$ contains $L_1 \bigcap (\bigcup_{i>1} L_i) \sub L_1$.  Then the conditions of \Cref{proposition:plumbings-immersions} are satisfied (take $f$ to be the identity map and let $K=L_1$). 
\epf

\ex[Plumbings]
Let $M_1,\dots, M_l$ be closed manifolds and let $X$ be the plumbing of the cotangent bundles of the $M_i$ according to any graph. Then $(X,\l)$ satisfies the assumptions of \Cref{proposition:plumbings}. In particular, $(X,\l)$ is big if $M_i \in \mathscr{E}$ for some $i$. 
\eex

In particular, a plumbing of cotangent bundles, one of which is (say) hyperbolic, is not symplectomorphic to an affine variety. 

Another rational natural class of skeleta arise from immersions. Here also, \Cref{proposition:plumbings-immersions} can be applied rather easily. For example:

\prop \label{proposition:immersions}
Suppose that $(X,\l)$ is Weinstein and that $\op{Skel}(X,\l)$ is the image of a Lagrangian immersion $f:\S_g \to (X,\l)$ with only simple double points, for $g \geq 2$. Then $(X,\l)$ is big. 
\eprop

\pf
Let $B \sub \S_g$ be a ball which contains the preimages of all the double points. Let $K$ be obtained by performing Polterovich surgery on the double points. Then $K$ is a (possibly non-orientable) Riemann surface which contains $\S_g-B$. By van Kampen, the image of the inclusion $\pi_1(\S_g- B) \to \pi_1(K)$ has exponential growth. Now apply \Cref{proposition:plumbings-immersions}. 
\epf

\rmk
The skeleton of the affine variety $\C^2- \{ xy=1\}$ is the \emph{Whitney immersion}, which is an immersion of the $2$-sphere with a single double point (see e.g.\ \cite{dr}). In contrast, it follows from \Cref{proposition:immersions} that the skeleton of an algebraic variety cannot be the image of an immersion of a surface of genus greater than $2$. We do not know whether there exists an algebraic variety whose skeleton is an immersed torus with at least one double point (of course, both $T^*S^2$ and $T^* \bb{T}^2$ are symplectomorphic to affine varieties). We remark that given any immersion $M \to \R^m$ with only simple double points, it is straightforward to construct a Weinstein manifold whose skeleton is the image of $f$.
\ermk

\subsection{A big Weinstein manifold having no exact Lagrangians}


For any $n \geq 6$, Abouzaid and Seidel \cite[Thm.\ 1.3]{abouzaid-seidel-recombination} constructed infinitely many examples of Weinstein manifolds of dimension $2n$ whose symplectic cohomology with $\Z/2$-coefficients vanishes but whose symplectic cohomology with $\Z$-coefficients (and hence also with $k$-coefficients) is nonzero. 

Let $X$ be such a Weinstein manifold and let $Y$ be a big Weinstein manifold. Then $X \tms Y$ is big by \Cref{proposition:products}. On the other hand, the K\"{u}nneth theorem for symplectic cohomology (see \cite{oancea}) implies that $SH^\bu(X \tms Y; \Z/2)= SH^\bu(X;\Z/2) \otimes SH^\bu(Y; \Z/2)=0.$ It follows by Viterbo restriction for symplectic cohomology with $\Z/2$-coefficients that $X \tms Y$ does not contain a closed exact Lagrangian.

\subsection{A big Weinstein manifold with boundary $S^5$} \label{subsection:exotics5}


Most of the arguments below were explained to the author by Marcelo Alves, and are drawn from \cite[Sec.\ 7.3]{alves-meiwes}. Consider the Brieskorn manifold $\S(2,3,7)$. This is an aspherical integral homology $3$-sphere; $\pi_1(\S(2,3,7))$ can be generated by two elements $\g_1, \g_2$, and this group has exponential growth (see \cite[Sec.\ 7.3]{alves-meiwes}).  

We consider $M:= S^* \S(2,3,7)= \S(2,3,7) \tms S^2$ (recall that all orientable $3$-manifolds are parallelizable). Let $Q= D^*\S(2,3,7)$. Choose isotropics representing $\g_1, \g_2$. 
Now attach handles along $\g_1, \g_2$. Let $Q'$ be the resulting Liouville domain, let $M'= \d Q'$ be the surgered contact manifold, and let $W$ be the trace of the surgery (i.e. $W$ is a cobordism with boundary $M' \sqcup \ov{M}$ and $Q'=Q \cup W$). By choosing the framing appropriately, we can assume that the first Chern class of $Q'$ is trivial (indeed, note that $c_1(Q)=0$ because $\S(2,3,7)$ is parallelizable; now the restriction of a trivialization of $\wedge_{\C}^3 TQ$ to the $\g_i$ determines the correct framing). 

\lem
We have $\pi_1(M')=0$ and $\pi_2(M')=H_2(M'; \Z) = \Z \oplus \Z$. 
\elem

\pf
The effect of surgery on homotopy and homology is well understood. To see that $\pi_1(M')=0$, simply apply \cite[Prop.\ 4.19(ii)]{ranicki} (put $m=5, n=1$ since we are dealing with  $1$-surgery on a $5$ manifold; hence the condition $2n+2 \leq m$ is verified). By Hurewicz, $\pi_2(M')=H_2(M'; \Z)$. Now $H_2(M'; \Z)$ can be computed by a routine diagram chase using \cite[Prop.\ 4.19(iii)]{ranicki} (again with $m=5, n=1$). 
%
%
\epf

By transversality, we can represent each element of $\pi_2(M')=H_2(M';\Z)= \Z^2$ by an embedded sphere. Let $\s_1, \s_2: S^2 \to M'$ be such embeddings which correspond to generators of $\Z^2$. 
By \cite[Lem.\ 2.6]{mclean-unpublished} and the fact that $c_1(Q')=0$, we can in fact assume that the $\s_1$ are (mutually disjoint) \emph{embedded Legendrians}.  Now let $\tilde{\s}_1, \tilde{\s}_2$ be the loosifications of $\s_1, \s_2$ (that is, the $\tilde{\s}_i$ are loose Legendrians which are embedded, disjoint, and $C^0$ close to the $\s_i$; their existence follows from Murphy's h-principle \cite[Thm.\ 7.25]{c-e}). 

Now attach handles along the $\tilde{\s}_i$; let $Q''$ be the resulting Liouville domain, let $M'' = \d Q''$ be the surgered contact boundary. One can verify that (use e.g.\  \cite[Lem.\ 1]{wall}) this handle attachment has the effect of killing $H_2(M';\Z)$; it also does not affect $\pi_1(M')=H_1(M';\Z)=0$ (use again e.g.\ \cite[Prop.\ 4.19(ii)]{ranicki}). Hence $\pi_1(M'')=H_1(M'';\Z)= H_2(M'';\Z)=0$. We conclude by appealing to the classical fact that a simply-connected closed $5$-manifold is entirely determined by its second homology (see \cite[Thm.\ 2.2]{barden}). In particular, it follows that $M''$ is the $5$-sphere. 

\cor
There exists a big Weinstein manifold whose ideal contact boundary is diffeomorphic to $S^5$.
\ecor
\pf
Indeed, $\widehat{Q''}$ is such a Weinstein manifold: it was constructed by first attaching subcritical handles to $Q$, and then attaching flexible handles.  Since $Q= D^*\S(2,3,7)$ and we saw that $\pi_1(\S(2,3,7))$ has exponential growth, it follows (\Cref{example:hyperbolics}) that $\widehat{Q}$ is big. But the property of being big is preserved under subcritical and flexible handle attachments (by \Cref{proposition:subcritical-embedding} and \Cref{proposition:flexible-faithful}), so the corollary follows. 
\epf

\cor
Any closed manifold of dimension $5$ which admits a Stein fillable contact structure admits a (possibly different) contact structure $\xi_{pos}$ with positive entropy (i.e. with the property that for every contact form inducing $\xi_{pos}$, the Reeb flow has positive topological entropy). 
\ecor
\pf
Let $V$ be closed contact manifold with Stein filling $W$. Then $W \sqcup \widehat{Q''}$ is certainly big. Now attach an appropriate (subcritical) $1$-handle whose effect on the contact boundary is to perform a connected sum of $V= \d_\infty W$ and $S^5= \d_\infty \hat{Q''}$. The resulting Weinstein manifold is still big, and $V \# S^5$ is of course diffeomorphic to $V$. 
\epf

\end{document}